\documentclass[11pt]{article}
\usepackage{geometry}                
\usepackage[parfill]{parskip}    
\usepackage{graphicx}
\usepackage{amsmath, latexsym, amssymb}
\usepackage{epstopdf}
\input xypic
\begin{document}
\title{Manifolds  admitting  stable forms}
\author{H\^ong-V\^an L\^e, Martin Panak  and Jiri  Vanzura}

\date{\today}                                           

\maketitle

\newcommand{\R}{{\mathbb R}}
\newcommand{\C}{{\mathbb C}}
\newcommand{\F}{{\mathbb F}}
\newcommand{\Z}{{\mathbb Z}}
\newcommand{\N}{{\mathbb N}}
\newcommand{\Q}{{\mathbb Q}}
\newcommand{\Hq}{{\mathbb H}}

\newcommand{\Aa}{{\mathcal A}}
\newcommand{\Bb}{{\mathcal B}}
\newcommand{\Cc}{{\mathcal C}}    
\newcommand{\Dd}{{\mathcal D}}
\newcommand{\Ee}{{\mathcal E}}
\newcommand{\Ff}{{\mathcal F}}
\newcommand{\Gg}{{\mathcal G}}    
\newcommand{\Hh}{{\mathcal H}}
\newcommand{\Kk}{{\mathcal K}}
\newcommand{\Ii}{{\mathcal I}}
\newcommand{\Jj}{{\mathcal J}}
\newcommand{\Ll}{{\mathcal L}}    
\newcommand{\Mm}{{\mathcal M}}    
\newcommand{\Nn}{{\mathcal N}}
\newcommand{\Oo}{{\mathcal O}}
\newcommand{\Pp}{{\mathcal P}}
\newcommand{\Qq}{{\mathcal Q}}
\newcommand{\Rr}{{\mathcal R}}
\newcommand{\Ss}{{\mathcal S}}
\newcommand{\Tt}{{\mathcal T}}
\newcommand{\Uu}{{\mathcal U}}
\newcommand{\Vv}{{\mathcal V}}
\newcommand{\Ww}{{\mathcal W}}
\newcommand{\Xx}{{\mathcal X}}
\newcommand{\Yy}{{\mathcal Y}}
\newcommand{\Zz}{{\mathcal Z}}

\newcommand{\zt}{{\tilde z}}
\newcommand{\xt}{{\tilde x}}
\newcommand{\Ht}{\widetilde{H}}
\newcommand{\ut}{{\tilde u}}
\newcommand{\Mt}{{\widetilde M}}
\newcommand{\Llt}{{\widetilde{\mathcal L}}}
\newcommand{\yt}{{\tilde y}}
\newcommand{\vt}{{\tilde v}}
\newcommand{\Ppt}{{\widetilde{\mathcal P}}}

\newcommand{\Remark}{{\it Remark}}
\newcommand{\Proof}{{\it Proof}}
\newcommand{\ad}{{\rm ad}}
\newcommand{\Om}{{\Omega}}
\newcommand{\om}{{\omega}}
\newcommand{\eps}{{\varepsilon}}
\newcommand{\Di}{{\rm Diff}}
\newcommand{\Pro}[1]{\noindent {\bf Proposition #1}}
\newcommand{\Thm}[1]{\noindent {\bf Theorem #1}}
\newcommand{\Lem}[1]{\noindent {\bf Lemma #1 }}
\newcommand{\An}[1]{\noindent {\bf Anmerkung #1}}
\newcommand{\Kor}[1]{\noindent {\bf Korollar #1}}
\newcommand{\Satz}[1]{\noindent {\bf Satz #1}}

\newcommand{\gl}{{\frak gl}}
\renewcommand{\o}{{\frak o}}
\newcommand{\so}{{\frak so}}
\renewcommand{\u}{{\frak u}}
\newcommand{\su}{{\frak su}}
\newcommand{\ssl}{{\frak sl}}
\newcommand{\ssp}{{\frak sp}}

\newcommand{\Cinf}{C^{\infty}}
\newcommand{\CS}{{\mathcal{CS}}}
\newcommand{\YM}{{\mathcal{YM}}}
\newcommand{\Jreg}{{\mathcal J}_{\rm reg}}
\newcommand{\Hreg}{{\mathcal H}_{\rm reg}}
\newcommand{\SP}{{\rm SP}}
\newcommand{\im}{{\rm im}}

\newcommand{\inner}[2]{\langle #1, #2\rangle}    
\newcommand{\Inner}[2]{#1\cdot#2}
\def\NABLA#1{{\mathop{\nabla\kern-.5ex\lower1ex\hbox{$#1$}}}}
\def\Nabla#1{\nabla\kern-.5ex{}_#1}

\newcommand{\half}{\scriptstyle\frac{1}{2}}
\newcommand{\p}{{\partial}}
\newcommand{\notsub}{\not\subset}
\newcommand{\iI}{{I}}               
\newcommand{\bI}{{\partial I}}      
\newcommand{\LRA}{\Longrightarrow}
\newcommand{\LLA}{\Longleftarrow}
\newcommand{\lra}{\longrightarrow}
\newcommand{\LLR}{\Longleftrightarrow}
\newcommand{\lla}{\longleftarrow}
\newcommand{\INTO}{\hookrightarrow}

\newcommand{\Sy}{\text{ Diff }_{\om}}
\newcommand{\Ex}{\text{Diff }_{ex}}
\newcommand{\jdef}[1]{{\bf #1}}
\newcommand{\QED}{\hfill$\Box$\medskip}

\newcommand{\UuU}{\Upsilon _{\delta}(H_0) \times \Uu _{\delta} (J_0)}
\newcommand{\bm}{\boldmath}

\begin{abstract}
 In this note we give a direct method to classify  all stable forms on $\R^n$  as well as to   determine their automorphism groups. We  show that  in dimension 6,7,8  stable forms coincide with non-degnerate forms.
 We present   necessary conditions and  sufficient conditions for a manifold to admit a stable form. We also discuss rich  properties  of the geometry of  such  manifolds.
\end{abstract}

\section{Introduction}

Special geometries defined by  a class of differential forms on manifolds are  again in  the center of interests of geometers.
These interests  are motivated by the fact that such a setting of special geometries   unifies  many known  geometries  as symplectic geometry and geometries 
with special holonomy [Joyce2000], as well as other geometries arised in the M-theory [GMPW2004],   [Tsimpis2005].  
A series of papers by Hitchin [Hitchin2000], [Hitchin2001]
and  his  school [Witt2005], etc., opened  a new way to these special geometries. Among them they studied geometries  associated with    certain stable 3-forms in dimensions 6, 7 and 8 
(see the definition of a stable form in section 2 after Proposition 2.2.)

To classify  the  stable forms on $\R^n$   one could use the classification  by Sato and Kimura  [S-K1977] 
of the  stable  forms  on $\C^n$  (they are  partial cases of
prehomogeneous spaces)  and to find the corresponding  real forms of the complex stable forms. It might  be possible  to define the automorphism groups of
stable forms on $\R^n$  by  using the  Sato-Kimura result for the complex case, but it might  not be straightforward, since in  general the structure of   algebraic ideals  over $\R$  is more complicated than the structure
of algebraic ideals  over $\C$. (The real problem here is to know how many connected components
the isotropy group of a  given 3-form  does have.)
We  also have noticed  a proof by Witt in  [Witt2005]  attempting to define   the automorphism group of
 the stable form of PSU(3)-type, but  unfortunately this proof is  incomplete (see  Remark 4.8 below).   
  
In sections 2, 3  we  study  some properties of  stable forms. In section 4 we classify  stable forms
on $\R^n$ and we determine their automorphism groups. Our classification is  based on  the Djokovic 
work [Djokovic1983].
In sections 5, 6, 7 we   present  certain necessary conditions as well  as some  sufficient conditions for a 
manifold to admit a stable form. We also discuss  the rich structure of manifolds admitting stable forms
in  sections 5, 6, 8. In particular we show that  for n = 7 or 8  the tangent bundle of any manifold $M^n$ which admits a stable 3-form
has a canonical structure of a  real simple Malcev algebra bundle.

\medskip

\section{Multi-symplectic forms and  stable forms}

We recall that a k-form  $ \gamma $  on  a vector space  $V^n$ over a field $F$ is called
{\bf multi-symplectic},  if  the following map
$$ I_{\gamma} :  V \to \Lambda ^{k-1} (V^n)^* : \,  v \mapsto  v \rfloor  \gamma$$
is  injective.

Clearly a 2-form is multi-symplectic, if and only if it is symplectic.

A  multi-symplectic form is generic in the following sense.
For any k-form $\gamma$ we can define its {\bf rank}, denoted by $\rho (\gamma)$, as the minimal
dimension of the subspace $ W \subset  V^*$ such that $\gamma  \in  \Lambda ^k W$.

\medskip

{\bf 2.1. Lemma.}  {\it   A k-form  $\gamma$ on  $ V ^n$ is  multi-symplectic, if and only if, its rank is
$n$.}

\medskip

{\it  Proof.}  It is easy to see that if  the rank of $\gamma $ is less than $n$ then
the linear map $I_\gamma$ has  non-trivial kernel.  On the other hand, if $I_\gamma$ has
 non-trivial kernel,  then $\gamma$ can be represented as a k-form in the dual space of the kernel.
 In fact we have that  the dimension of kernel of $I_\gamma$  is equal to $n - \rho (\gamma)$\QED
 
 From now on
we shall assume that $F = \C$ or  $\R$.  In these cases the space $\Lambda^k  ( V^n)^* $ has
the natural topology induced from   $F$.  

\medskip

{\bf 2.2. Proposition}. {\it  The set of  multi-symplectic  k-forms  is open and dense in the
space of all k-forms.}

\medskip

{\it Proof.}  The equation for $\gamma \in  \Lambda^k(V^n) ^*$ defining that  $I_\gamma$ has  non-trivial  kernel  is an 
algebraic equation, so  the set of non-multi-symplectic k-forms is  a closed subset in $\Lambda
 ^k (V^n) ^*$.  It is  also easy to check that for any $k$ there exists a multi-symplectic k-form
 on $V^n$. Hence follows the statement.\QED

Clearly  the multi-symplecity is invariant under the action of the group $GL (F^n)$.  
We shall say that a k-form  $\gamma$ is {\bf stable}, if the orbit $GL(F^n)(\gamma)$
is open in the space  $\Lambda ^k (V^n) ^*$.  By Proposition 2.2 the set of multi-symplectic k-forms
has  non-trivial intersection with the orbit of any stable form. Hence follows immediately

\medskip

{\bf 2.3.  Corollary.}  {\it  A stable form is multi-symplectic.}

\medskip

  The  converse statement is true for  $k = 2$ or $k = n -2$.  If $k=3$ and
 $n = 7$, $ F = \R$, it  is known  that  there are  8 types of   $GL (\R ^7)$-orbits of multi-symplectic 3-forms but
 among them there are only two of them are stable.

\medskip

We say that two forms are equivalent (or of the same type),  if they are in the same orbit of $GL (V^n)$-action.  Clearly  a real form is stable, if and only if its complexification
is  stable. We also know that each complex orbit has a finite number of real forms [B-C1962, Proposition 2.3].
Thus the classification of  real stable forms is  equivalent to the classification of 
complex stable forms  plus the classifacation of the real forms of the complex stable forms. 
The classifications of  complex stable forms is a part of the  Sato-Kimura classification of
prehomogeneous spaces [S-K1977].  

\medskip

\section {Symmetric  bilinear forms associated  to  a 3-form on $\R^8$.}

In this section we associate  to  a
3-form  $\om^3$ on $\R^8$ several  symmetric bilinear forms which are invariants of  $\om^3$.  
We  prove that the only
 non-degenerate 3-forms (see  definition  below, after formula (3.4))  are stable forms. For each stable form we shall associate a Lie
 algebra structure on $\R^8$.

\medskip

We denote by $I$ the  following natural isomorphism  from $  \R^8 \otimes \Lambda ^8 (\R ^8)^* \to  \Lambda ^7 ( \R^8)^*$:
$$ I ( v\otimes  \theta )  = v \rfloor  \theta. \leqno(3.1)$$

Let $\om $ be a 3-form on $\R^8$. We associate $\om$ with a  symmetric
bilinear map $S : \R ^8 \times  \R^8 \to  \R^8 \otimes \Lambda ^8 (\R^8)^*$ as follows
$$ S^\om( v, w) =  I ^{-1} ((v \rfloor \om) \wedge  (w\rfloor \om) \wedge \om ). \leqno (3.2) $$
Equivalently
$$ S^\om (v,w) =  - \sum_{i=1}^8 e _i \otimes ((v \rfloor \om) \wedge (w\rfloor \om) \wedge \om  \wedge e_i ^*)\leqno (3.2.a)$$
 for any basis $(e_i)$ in  $ \R^8$ and its dual basis $(e_i ^*)$. 
 
For each $v \in \R^8$ we  define a linear map $ L ^\om _v : \R^8 \to \R^8  \otimes \Lambda ^8  (\R^8)^*$ by  letting the first   variable in $S^\om$ to be $v$
$$ L^\om _v (w) = S^\om (v, w). \leqno (3.3)$$
Now we  shall define a symmetric linear  form $B^\om (v,w): \R^8 \times \R^8 \to ( \Lambda ^8 (\R^8) ^*) ^2$  as follows
$$B^\om (v,w) =  Tr (L^\om _v  \circ L ^\om_w) \in  ( \Lambda ^8 (\R^8) ^*) ^2. \leqno (3.4)$$
We  say that $\om$ is {\bf non-degenerate}, if  the reduced trace form $<B^\om, \rho^2>$ is non-degenerate, for some  choice  of $\rho \in \Lambda ^8 (\R^8) \setminus \{ 0\}$.

Let $G_\om$ be the automorphism group  of $\om$.  
Let us  consider  the  component   $G_\om ^+: =G_\om  \cap Gl ^+ ( \R^8)$.

\medskip

{\bf 3.5. Proposition}.  {\it  The   bilinear  forms  $S^\om$ and  $B^\om$ are   $Gl (\R^8)$-equivariant
  in the following sense. For any $g \in  Gl (\R^8)$ we have
$$  S^{g ^*(\om)}(X, Y) =  g^* ( S^{\om} (g^{-1} X,  g^{-1} Y) ),\leqno (3.5.1) $$
$$B^{g^* (\om)} (X, Y) = g  ^* ( B^\om (g^{-1}  X, g^{-1} Y)). \leqno (3.5.2)$$
If  $\om$ is non-degenerate, then the group $G_{\om} ^+$ is  a subgroup of $ SL (\R^8)$. The group $G_\om$ preserves  the reduced  trace form  $<B^\om  , \rho^2>$  for any choice of  $\rho \in \Lambda ^8( \R^8)$. }

\medskip

{ \it Proof}.  The computation of (3.5.1) and (3.5.2) is straightforward, so we omit them.
 The symmetric form $B^\om(v,w)$ can be considered as a linear map
$B^\om : ( \R^8 )   \to  (\R^8)^* \otimes (\Lambda ^8 (\R^8)^*) ^2$. Let us consider the associated linear map
$$\det  ( B^\om):   \Lambda  ^8 ( \R^8) \to \Lambda^8((\R^8)^* \otimes (\Lambda ^8 ( \R^8) ^*) ^2) = \Lambda ^8( (\R^8)^*) ^{17}.\leqno(3.5.3)$$
If  $B^\om$ is non-degenerate, then the map $\det (B^\om)$ is not trivial.   From (3.5.2) we  deduce that  the map $\det B^\om$ is $G_\om ^+$-invariant map.  So for any  $g \in G_\om ^+$ we get from (3.5.3)
$$ \det  g = (\det g ^{-1} ) ^{17}.$$
Since $\det g  > 0$ we conclude that  $\det g = 1$.  Now using (3.5.2) we get the last   statement immediately.\QED

\medskip

{\bf 3.6. Proposition.}  {\it i) The trace form $B^\om$  is compatible  with the multiplication
$S^\om$  in the following sense
$$B^\om( S^\om (a.b), c) =  B^\om ( a, S^\om (b,c)).$$
ii) The trace form $B^ \om$ is non-degenerate, if and only if $\om$ is  stable.}

\medskip

{\it Proof.}   The first statement follows immediately from the definition. To prove the second statement
we observe that if $\om_1$ and $\om_2$  are the real forms of the same complex 3-form,
then  their trace forms are also  the real forms of the trace form for the complex 3-form (all these
bilinear forms $S^\om$ and $B^\om$ can be defined for any vector space $V$ over an arbitrary field.)
Thus to check  how many real 3-forms are non-degenerate we need to check only 22 representatives
of 3-forms in the Djokovic classification [Djokovic1983].  Furthermore  we know that  a non-degenerate
3-form must be multi-symplectic. Thus it suffices to compute the trace form of  13 
multi-symplectic 3-forms
 in tables XI-XXIII in  the Djokovic classification. We wrote a program for computing  the trace
form  $B^\om$  to run it under Maple.   The program is very simple. 
We denote by $e_1^*\wedge \cdots\wedge e_8^*$ by  $\theta$, where
$e_i^*$ are the coordinate 1-forms  on $\R^8$.
We shall use  $\theta$ to make a (reduced) multiplication $V \times V \to V$
$$(vw  \rfloor \theta) = (v\rfloor \om) \wedge (w \rfloor \om) \wedge \om \leqno (3.7) $$
Clearly we have
$$ S^\om (v, w) =  vw \otimes \theta.\leqno (3.8)$$

We define structure constants $A_{ij}^k$ by
$$ e_i e_j = \sum _k  A_{ij}^k e_k \leqno (3.9)$$
Then 
$$ S^\om ( e_i, e_j)  = \sum_k A_{ij}^k e_k \otimes  \theta \leqno (3.9.a)$$
Now let us compute
$$B^\om ( e_l, e_m) =  \sum_n  ( S (e_l, S (e_m, e_n)), e_n ^*) $$
$$\stackrel {3.2.a}{ = } \sum_{k ,n}  <e_k \otimes  ( e_l \rfloor \om)\wedge( e_m e_n \rfloor \om)  \wedge \om
\wedge e_k ^* \otimes \theta, e_n ^*> $$
$$ =  \sum_{n, p} ( e_l \rfloor  \om) \wedge A_{mn} ^p  (e_p \rfloor  \om) \wedge \om \wedge e_n^*\otimes \theta$$
$$ = \sum _{n, p} A _{lp}^n \cdot A_{ m,n} ^p \otimes (\theta)^2. \leqno(3.10)$$

The result is that the only stable forms numerated by XXIIIa, XXIIIb,
XXIIIc by Djokovic  have non-degenerate  trace forms.

\medskip

Below we shall compute explicitly   the reduced multiplication  forms as well as 
the reduced trace  forms $<B^{\phi_i}, (\theta ^*) ^2>$ for  stable  forms 
$\phi_i$ on $\R^8$
from the Djokovic classification.

(Form XXIIIa): $\phi_1 = e^{124} + e^{134}
 +e^{256} + e ^{378} + e ^{157} + e ^{468}$.\\
 (Form XXIIIb): $\phi_2 =
 e^{135} + e ^{245} + e ^{146}  - e ^{236} + e ^{127} + e ^{348} + e ^{678}$.\\
(Form XXIIIc): $\phi _3 = e ^{135}  -e ^{146} + e ^{236} + e ^{245} + e ^{347} + e^{568} + e ^{127}
 + e ^{128}$. 
 
The  reduced multiplication table for  the form XXIIIa is: 
$$
\left[ \begin {array}{cccccccc} 0&-{\it e1}&{\it e1}&3\,{\it e2}-3\,{\it e3}&-3\,{\it e8}&0&-3\,{\it e6}&0\\\noalign{\medskip}-{\it e1}&-2
\,{\it e2}&-2\,{\it e2}+2\,{\it e3}&-{\it e4}&-{\it e5}&-{\it e6}&2\,{
\it e7}&2\,{\it e8}\\\noalign{\medskip}{\it e1}&-2\,{\it e2}+2\,{\it 
e3}&2\,{\it e3}&{\it e4}&-2\,{\it e5}&-2\,{\it e6}&{\it e7}&{\it e8}
\\\noalign{\medskip}3\,{\it e2}-3\,{\it e3}&-{\it e4}&{\it e4}&0&0&3\,
{\it e7}&0&3\,{\it e5}\\\noalign{\medskip}-3\,{\it e8}&-{\it e5}&-2\,{
\it e5}&0&0&3\,{\it e3}&-3\,{\it e4}&0\\\noalign{\medskip}0&-{\it e6}&
-2\,{\it e6}&3\,{\it e7}&3\,{\it e3}&0&0&3\,{\it e1}
\\\noalign{\medskip}-3\,{\it e6}&2\,{\it e7}&{\it e7}&0&-3\,{\it e4}&0
&0&-3\,{\it e2}\\\noalign{\medskip}0&2\,{\it e8}&{\it e8}&3\,{\it e5}&0
&3\,{\it e1}&-3\,{\it e2}&0\end {array} \right] 
$$

The reduced trace form for the form XXIIIa  is:
$$
\left[ \begin {array}{cccccccc} 0&0&0&-30&0&0&0&0\\\noalign{\medskip}0&20&10&0&0&0&0&0\\\noalign{\medskip}0&10&20&0&0&0&0&0
\\\noalign{\medskip}-30&0&0&0&0&0&0&0\\\noalign{\medskip}0&0&0&0&0&-30
&0&0\\\noalign{\medskip}0&0&0&0&-30&0&0&0\\\noalign{\medskip}0&0&0&0&0
&0&0&-30\\\noalign{\medskip}0&0&0&0&0&0&-30&0\end {array} \right] 
$$

The reduced multiplication table for the form XXIIIb is: 
$$
 \left[ \begin {array}{cccccccc} 6\,{\it e8}&0&-3\,{\it e6}&3\,{\it e5
}&-{\it e1}&3\,{\it e2}&-3\,{\it e4}&0\\\noalign{\medskip}0&6\,{\it e8
}&-3\,{\it e5}&-3\,{\it e6}&-{\it e2}&-3\,{\it e1}&3\,{\it e3}&0
\\\noalign{\medskip}-3\,{\it e6}&-3\,{\it e5}&6\,{\it e7}&0&-{\it e3}&
-3\,{\it e4}&0&3\,{\it e2}\\\noalign{\medskip}3\,{\it e5}&-3\,{\it e6}
&0&6\,{\it e7}&-{\it e4}&3\,{\it e3}&0&-3\,{\it e1}
\\\noalign{\medskip}-{\it e1}&-{\it e2}&-{\it e3}&-{\it e4}&-2\,{\it 
e5}&2\,{\it e6}&2\,{\it e7}&2\,{\it e8}\\\noalign{\medskip}3\,{\it e2}
&-3\,{\it e1}&-3\,{\it e4}&3\,{\it e3}&2\,{\it e6}&-6\,{\it e5}&0&0
\\\noalign{\medskip}-3\,{\it e4}&3\,{\it e3}&0&0&2\,{\it e7}&0&0&3\,{
\it e5}\\\noalign{\medskip}0&0&3\,{\it e2}&-3\,{\it e1}&2\,{\it e8}&0&
3\,{\it e5}&0\end {array} \right] 
$$

The reduced trace form for the form XXIIIb is:
$$
\left[ \begin {array}{cccccccc} 0&0&0&-60&0&0&0&0\\\noalign{\medskip}0&0&60&0&0&0&0&0\\\noalign{\medskip}0&60&0&0&0&0&0&0
\\\noalign{\medskip}-60&0&0&0&0&0&0&0\\\noalign{\medskip}0&0&0&0&20&0&0
&0\\\noalign{\medskip}0&0&0&0&0&-60&0&0\\\noalign{\medskip}0&0&0&0&0&0
&0&30\\\noalign{\medskip}0&0&0&0&0&0&30&0\end {array} \right] 
$$

The  reduced multiplication  table for the form XXIIIc is: 

$$
\left[ \begin {array}{cccccccc} 6\,{\it e7}-6\,{\it e8}&0&3\,{\it e6}&3\,{\it e5}&3\,{\it e4}&3\,{\it e3}&{\it e1}&-{\it e1}
\\\noalign{\medskip}0&6\,{\it e7}-6\,{\it e8}&-3\,{\it e5}&3\,{\it e6}
&-3\,{\it e3}&3\,{\it e4}&{\it e2}&-{\it e2}\\\noalign{\medskip}3\,{
\it e6}&-3\,{\it e5}&6\,{\it e8}&0&-3\,{\it e2}&3\,{\it e1}&{\it e3}&2
\,{\it e3}\\\noalign{\medskip}3\,{\it e5}&3\,{\it e6}&0&6\,{\it e8}&3
\,{\it e1}&3\,{\it e2}&{\it e4}&2\,{\it e4}\\\noalign{\medskip}3\,{
\it e4}&-3\,{\it e3}&-3\,{\it e2}&3\,{\it e1}&-6\,{\it e7}&0&-2\,{\it 
e5}&-{\it e5}\\\noalign{\medskip}3\,{\it e3}&3\,{\it e4}&3\,{\it e1}&3
\,{\it e2}&0&-6\,{\it e7}&-2\,{\it e6}&-{\it e6}\\\noalign{\medskip}{
\it e1}&{\it e2}&{\it e3}&{\it e4}&-2\,{\it e5}&-2\,{\it e6}&2\,{\it 
e7}&2\,{\it e7}-2\,{\it e8}\\\noalign{\medskip}-{\it e1}&-{\it e2}&2\,
{\it e3}&2\,{\it e4}&-{\it e5}&-{\it e6}&2\,{\it e7}-2\,{\it e8}&-2\,{
\it e8}\end {array} \right] 
$$

The reduced trace 
form for the form XXIIIc is:
$$
\left[ \begin {array}{cccccccc} 60&0&0&0&0&0&0&0\\\noalign{\medskip}0&60&0&0&0&0&0&0\\\noalign{\medskip}0&0&60&0&0&0&0&0
\\\noalign{\medskip}0&0&0&60&0&0&0&0\\\noalign{\medskip}0&0&0&0&60&0&0
&0\\\noalign{\medskip}0&0&0&0&0&60&0&0\\\noalign{\medskip}0&0&0&0&0&0&
20&10\\\noalign{\medskip}0&0&0&0&0&0&10&20\end {array} \right] 
$$
\QED

\medskip

{\bf 3.11.  Proposition.} {\it Each  stable form $\phi$ defines a Lie algebra structure $[,]_\phi$ on  $\R^8$ by the 
following formula

$$ <[X, Y]_\phi ,Z>_\phi  = \phi (X, Y, Z ),\leqno (3.11.1)$$
where $<, >_\phi$ denotes a  reduced trace form of $\phi$.
Moreover the Lie algebra  $[,]_{\phi_i}$  is the non-compact  real form of
$sl(3,\C)$ for $i =1,2$ and  the Lie algebra  $[,]_{\phi_3}$ is the compact real form
of $sl (3, \C)$.}

\medskip

{\it Proof.} First we  note that the anti-symmetric bracket $[,]_\phi$ satisfies the following invariant property.
For each $g \in Gl (\R^8)$ we have
$$  [X, Y]_{g ^* {\phi}} = g ([g^{-1}(X), g^{-1} (Y)]) _\phi.\leqno (3.12)$$
Hence if the Jacobi identity holds at a form $\phi$, it also holds at any point in the orbit
$GL(\R^8) (\phi)$, moreover these Lie brackets are equivalent.
Secondly we notice that the bracket $[,]_\phi$  can be extended  linearly over $\C$ and this complexification is
the anti-symmetric bracket defined by the complexification of the form $\phi$ according to the same formula  (3.11.1). 
Thus to verify the Jacobi identity  for  3 stable forms $\phi_i, i =1,3$, it suffices to verify for one of them.

Now we have two proofs  for Proposition 3.11. In the first way we   compute the Lie bracket defined in (3.11.1)  
 by using  our   explicit formula for the reduced
trace forms of  one of  the real stable form $\phi$ and   arrive at the above conclusion. 
The second method  uses the Cartan form on
the real form  $su(3)$ of  the  complex Lie algebra $sl (3, \C)$. 
First 
we  compute
the reduced trace formula for the Cartan form on the algebra $su(3)$
$$\rho _3 (X, Y, Z) = <[X,Y], Z>$$
where $<,>$ denotes the Killing form on $su (3)$.  We use the following explicit expression taken from [Witt2005] for 
a multiple
of the form $\rho_3$:
$$ (-1/\sqrt 3)^3\rho_3 = e^{123}+(1/2) (e^{147}-e^{156}+ e^{246} + e ^{257} + e^{345}-e^{367}) + (\sqrt 3/2) 
(e^{845} + e^{867})$$
where $(e_i)$ are an orthonormal  basis in $su (3)$ and $e^{ijk}$ denotes the form $e^i \wedge e^j \wedge e^k$.
A direct computation (also used  Maple) gives us the following multiplication table 
for  $(4/3) \cdot (-1/\sqrt 3)^3\rho_3 $
$$
\left[ \begin {array}{cccccccc} 2\,{\it e8}&0&0&\,\sqrt {3}{\it e6}&\,\sqrt {3}{\it e7}&\,\sqrt {3}{\it e4}&\,\sqrt {3}{\it 
e5}&2\,{\it e1}\\\noalign{\medskip}0&2\,{\it e8}&0&-\,\sqrt {3}
{\it e7}&\,\sqrt {3}{\it e6}&\,\sqrt {3}{\it e5}&-\,\sqrt {3}
{\it e4}&2\,{\it e2}\\\noalign{\medskip}0&0&2\,{\it e8}&\,
\sqrt {3}{\it e4}&\,\sqrt {3}{\it e5}&-\,\sqrt {3}{\it e6}&-
\,\sqrt {3}{\it e7}&2\,{\it e3}\\\noalign{\medskip}\,\sqrt {3}{
\it e6}&-\,\sqrt {3}{\it e7}&\,\sqrt {3}{\it e4}&\,\sqrt {3}{
\it e3}-\,{\it e8}&0&\,\sqrt {3}{\it e1}&-\,\sqrt {3}{\it e2}
&-\,{\it e4}\\\noalign{\medskip}\,\sqrt {3}{\it e7}&\,\sqrt {
3}{\it e6}&\,\sqrt {3}{\it e5}&0&\,\sqrt {3}{\it e3}-\,{\it 
e8}&\,\sqrt {3}{\it e2}&\,\sqrt {3}{\it e1}&-\,{\it e5}
\\\noalign{\medskip}\,\sqrt {3}{\it e4}&\,\sqrt {3}{\it e5}&-
\,\sqrt {3}{\it e6}&\,\sqrt {3}{\it e1}&\,\sqrt {3}{\it e2}&-
\,\sqrt {3}{\it e3}-\,{\it e8}&0&-\,{\it e6}\\\noalign{\medskip}
\,\sqrt {3}{\it e5}&-\,\sqrt {3}{\it e4}&-\,\sqrt {3}{\it e7}
&-\,\sqrt {3}{\it e2}&\,\sqrt {3}{\it e1}&0&-\,\sqrt {3}{\it 
e3}-\,{\it e8}&-\,{\it e7}\\\noalign{\medskip}2\,{\it e1}&2
\,{\it e2}&2\,{\it e3}&-\,{\it e4}&-\,{\it e5}&-\,{\it e6}&
-\,{\it e7}&-2\,{\it e8}\end {array} \right] 
 $$

and  we compute easily from here (also  by using  Maple) that
 the reduced  trace formula for
$(-1/\sqrt 3)^3\rho_3 $ is equal to $(45/4)\,  (diag)$.  So the trace formula is a multiple  of the Killing form.
In particular it is non-degenerate. From Proposition 3.6 we obtain that the  Cartan form  on $su(3)$ is stable.
Now we observe that  the other Cartan forms   on $sl(3, \R)$ and on $su (1,2)$ are the real forms of the Cartan form
on $sl (3,\C)$ so they are also stable. 

(There is also another argument to prove that the Cartan forms are stable without using  Proposition 3.6 and using the
Djokovic classification. We need only to compute
their trace form explicitly and notice that the reduced trace form is a multiple of  the Killing form. Then we apply  
the argument
we use in the proof of the classification theorem 4.1 to get the dimension of the isotropy group of  the Cartan form and 
hence we get  the stability of the Cartan form.)

Once we know that the reduced trace form is a multiple of the Killing form, we get immediately   the first statement of Proposition 3.11.  

Clearly  the Lie algebra $g$ must lie in the Lie algebra of the stabilzer
of the Cartan form. Comparing  these Lie algebras with the Lie algebra of the automorphism group of the stable forms in  the Djokovic classification
we conclude that the Cartan form on $su (3)$ is equivalent with the form $\phi_3$, the Cartan form on $sl(3,\R)$
is equivalent to the form $\phi_1$ and the Cartan form on $su(1,2)$ is equivalent to the form $\phi_2$.
This proves the second statement of Proposition 3.11.\QED

\medskip

\section {Classification of real stable forms.}

We observe that  the stability of a k-form is preserved  under the Hodge isomorphism
$\Lambda ^k ( V^n) ^* \to \Lambda  ^{n-k} ( V^n) ^*$. We shall   use notation
$ e^{12\cdots k}$ for the form  $e^1 \wedge e ^2 \wedge \cdots \wedge e^k$. 
We also use notation  $G_\gamma$ for  the isotropy group of  $\gamma$ under
the action of $Gl (\R^n)$ and by $g_\gamma$ the Lie algebra of $G_\gamma$.

\medskip

{\bf 4.1. Theorem}. {\it  Suppose that $3 \le  k \le  n-k$. \\
i)Then a stable  k-form $\gamma$ on
 $\R ^n$   exists,  if and only if $k =3$ and $6 \le n\le 8$.  Furthermore  \\
ii) if  n = 6, then $\gamma$ is equivalent to one of the following forms:\\
 $\gamma_1 =e^1 
 \wedge e^2 \wedge e^3 + e^4\wedge e^5 \wedge e^6$  with
 $G_{\gamma_1} = ( SL(\R^3) \times SL (\R^3) \times ) \Z_2$;\\
 $\gamma_2 =  Re \, (e^1 +i e^2 ) \wedge (e^3 + i e^4) \wedge (e^5 + i e^6)$
 with $G_{\gamma _2} = SL(\C ^3)$,\\
 iii) if n = 7, then  $\gamma$ is equivalent to  one of the following forms: \\
  $ \om_1 = e ^{123} - e^{145} + e^{167}  + e^{246} + e^{257} + e ^{347} - e^{356}$ with $G_{\om _1} = G_2$;\\
    $ \om_2 = e ^{123} + e^{145} - e ^{167} +
 e^{246} + e^{257} + e^{347} - e ^{356}$ with $G_{\om_2} = \tilde G _2$,\\
 iv)  if  n =8, then $\gamma $ is equivalent to one of the following forms: \\
  $\phi_1 = e^{124} + e^{134}
 +e^{256} + e ^{378} + e ^{157} + e ^{468} $ with  $G_{\phi_1} = SL (3, \R)\times) \Z_2$;\\
  $\phi_2 =
 e^{135} + e ^{245} + e ^{146}  - e ^{236} + e ^{127} + e ^{348} + e ^{678}$ with  $g_{\phi_2}
 = PSU (1,2)\times) \Z_2$;\\
   $\phi _3 = e ^{135}  -e ^{146} + e ^{236} + e ^{245} + e ^{347} + e^{568} + e ^{127}
 + e ^{128}$ with $G_{\phi_3}  = PSU(3) \times) \Z_2$. }
 
\medskip

{\it Proof.} We first show that if $4 \le k \le n -k$ then  there is no stable form.  
It suffices to show  that
 in this case
 we have
 $$ \dim \Lambda ^k  (  \R^n) \ge n ^2 +1 = \dim (Gl ( V^n)) + 1. \leqno (4.2)$$
 Clearly we have under the  assumption that $4 \le k \le n-k$
 $$\dim \Lambda ^k (\R^n)  \ge  \dim \Lambda ^4 (\R^n).$$
 Therefore (2.2) is a consequence of the following equality
 $$f(n) : =  n^3 -6n ^2 -13 n - 6  \ge 1,  \text  { for  } n  \ge 8. \leqno (4.3)$$
Since $f' (n) > 0$ for all $n \ge 8$  it suffices to check (4.3) for $n = 8$ which is an easy exercise.
To complete  the proof of Theorem 4.1.i we need to show that stable 3-forms exist
for $n =6,7,8$ and not for $n \ge 9$.  
 But this is  an well-known fact for $n = 6,7$ and it follows from the classification
of 3-forms on  $\R^8$ by Djokovic  [Djokovic1983].  To show that there is no stable
3-form in $\R^n$, if $n \ge 9$ we can repeat the argument above to show that  in this case
$\dim \Lambda ^3 (\R^n)  > \dim Gl(\R^n)$.

ii) This classification is already well-known,  see  [Hitchin2000] for an wonderful treatment.

iii) This classification  follows from  the list of Bures and Vanzura of multi-symplectic
3-forms in dimension 7 [B-V2003] together with their  automorphism groups.  The  groups $G_{\om_i}$  have been first determined by
Bryant   [Bryant1987].

iv) We shall  complete this classification   from  the last table in [Djokovic1983]. In that table
Djokovic  supplied us only the Lie algbras $g_{\phi_i}$, for $i =1,2,3$.  We shall recover  $G_{\phi_i}$ from
$g_{\phi_i}$ by  using the following lemmas 4.4 and 4.5.  

\medskip

{\bf 4.4. Lemma}.   {\it Group $Gl^+( \R^n)$ acts transitively on the orbit 
$Gl (\R^n) ( \phi_i)$,  for  $\phi_i$ being one of the forms in Theorem 2.1.iv.}

\medskip

{\it Proof.}  It suffices to show that  the intersection $G_{\phi_i}\cap Gl^-(\R^8)$ is not empty, 
where $Gl^-(\R^8)$ denotes the orientation reserving  component of $Gl (\R^n)$. 

- For $\phi_1$ this intersection contains the following element  $\sigma_{23} \cdot \sigma_{57}\cdot
\sigma_{68}\cdot I_1 \cdot I _4$. Here  $\sigma_{ij}$ denotes the orientation reserving linear transformation which permutes the basic vectors $v_i$ and $v_j$ and leaves all other  basic vectors,
and $I_j$ denotes the  orientation reserving linear transformation which acts as $-Id$ on the
 line $v_j \otimes \R$ and leaves all other basic vectors. 
 
 - For  $\phi_2$ this intersection contains the following element  $\sigma_{12}\cdot \sigma_{34}\cdot
 I_6 \cdot I_7 \cdot I_8$.
 
 - For   $\phi_3$ this intersection contains $\sigma_{34}\cdot \sigma_{56}\cdot I_1\cdot I_7 \cdot I_8$.
 \QED

{\bf 4.5. Lemma.}  {\it Group $Gl_{\phi}^+  : = Gl^+ (\R^8) \cap G_{\phi}$ is connected  for  $\phi_i$
being one of the forms in Theorem 4.1.iv.}

\medskip

{\it Proof.}  We use the observation obtained in  section 3 that
 all three forms $\phi_i$  are the Cartan forms
 $$ \rho(X, Y, Z) =  < X, [Y, Z ] >$$
 on  the Lie algebra $sl (3, \R)$, $su (1,2)$ and $su (3)$, where $<,>$ denotes the Killing form. 
 Hence follows  that 
 $$Aut (g_{\phi_i}) \subset G_{\phi_i}.\leqno (4.6)$$
 
In Proposition 3.11 we  have also defined a way to recover the  structure of the corresponding Lie algebra  from $\phi_i$.  Since  all the reduced  bilinear forms are invariant  with respect to
$G_{\phi_i } $ we get
$$G_{\phi_i} \subset  Aut (g_{\phi_i}).\leqno(4.7)$$
Finally  the structure of $Aut (g_{\phi_i})$ is well-known, see e.g. [Murakami1952] and the references therein.  Thus we get Lemma 4.5 from  (4.6) and (4.7).\QED

Actually  the proof of Lemma 4.5 implies  Lemma 4.4.  Nevertheless the proof of Lemma
4.4 gives us explicitly an element in $Gl^- (\R^8) \cap G_{\phi_i}$. 
This completes the proof of Theorem 4.1.\QED

\medskip

{\bf 4.8. Remark}. In his thesis [Witt2005]  Witt  gave a proof that the component $G_{\phi_3}^+$
is $PSU(3)$.  His proof is incomplete, since he used implicitly without a proof
the fact  that  the component $G_{\phi_3}^+$  preserves the Killing metric on $su (3)$.  (His method
is to associate   the  Cartan form  to a  bilinear form  with value on $\R ^8$ by using a fixed   basis
of  $\R^8$.  A detailed analysis shows that such a  use is  equivalent to giving a linear  map from
$(\R^8) ^*$ to  $\R^8$ and in the  given case of Witt, that map is an isomorphism defined by the Killing metric).

\medskip

We say that   a differentiable form $\gamma$ on  a manifold $M^n$ is stable, if
at each $x \in M$ the form $\gamma(x)$  is stable.

\medskip

{\bf 4.9. Proposition.}  {\it   If a connected  manifold  $M^n$ admits a differentiable stable   form  $\gamma^3$,
then  for all $x\in M^n$ the form  $\gamma(x)$ has the same type. In particular $M^n$ admits
a $G_{\gamma (x) }$  structure. Conversely, if  $M^n$ admits a  $G_{\gamma}$ structure, then
it admits a differentiable form of $\gamma$ type.}

\medskip

{\it Proof.}  For each $x\in M^n$ denote by $U(x)$ the set of all  points $y \in M^n$ such that
$\gamma^3(y)$ has the same type as $ \gamma ^3(x)$. Clearly $ U(x)$ is   an open
subset in $M^n$.  Suppose that $U(x)\not = M^n$. Then  the  closure  $\bar U(x)$  contains an point $y$ which is not
in $U(x)$.  Clearly $\gamma(y)$ also has the same type as $\gamma(x)$ since
$U(y)$ has a non-empty intersection with $U(x)$. Thus $y \in U(x)$  which is a contradiction. The last statement  follows from the fact that  the transition functions on $G(x)$-manifold preserve
the form $\gamma (x)$.\QED

\medskip

\section{Stable 3-forms on 6-manifolds}

{\it 5.1. Obstruction for the existence of a $G_{\gamma_1}$-structure.}

If a non-orientable manifold $M^6$ admits a  $G_{\gamma_1}$-structure, then its orientable double covering
shall admit  $G_{\gamma_1}$-structure.  Now we shall concern only orientable manifolds $M^6$ and so only the identity component  of $G_{\gamma_1}$.
Clearly  $M^6$ admits $SL(3)\times SL(3)$-structure, if
and only if it admits a distribution of  oriented 3-planes on $M^6$.

We denote by $\rho_2 : H^2(M, \Z) \to  H^2(M, \Z_2)$ the modulo 2 reduction.

\medskip

{\bf 5.1.1. Proposition}.  {\it  Suppose that a closed manifold $M^6$  admits a $SL(3)\times SL(3)$-structure. Then its Euler class vanishes. Assume  that   $H^4(M^6, \Z)$ has no 2-torsion, the Euler class
$e(M^6)$ vanishes and
$M^6$  satisfies moreover the following condition (P). There are classes $c_1, c_2 \in H^2(M, \Z)$
such that
 $$ p_1 (M^6) = c_1^2 + c_2^2, \:  \rho_2 (c_1 + c_2) = w_2 (M^6). \leqno (P)$$
Then $M^6$ admits a $SL(3)\times SL(3)$-structure.}
\medskip

{\it Proof}.   The first statement is well-known, since  the  Euler class of an oriented 3-dimensional vector bundle is a 2-torsion, and $H^6(M, \Z)$ has no 2-torsion.
Let us  assume that  an  orientable manifold  $M^6$ with vanishing Euler class has no 2-torsion in $H^4(M, \Z)$, moreover $M^6$ satisfies condition (P).  Let $V$ be a non-vanishing vector field on $M^6$. 
Since $M^6$ satisfies condition (P), there is an almost complex structure $J$ on $M^6$ such that
$c_1(J) = c_1 +c_2$, where $c_1$ and $c_2$ satisfies condition (P).  Let $W^4$ be a   $J$-invariant   sub-bundle of $TM^6$ which
is complement to $V$ and $JV$.   Clearly $ p _1(W^4) = p_1 (M^6)$. Let $L_1$ and $L_2$ be
the complex line bundles over $M^6$  with the first  Chern classes $c_1$ and $c_2$ satisfying condition
(P).
Then $p_1(W^4) = p_1 (M^6) = p_1 (L_1 \oplus L_2)$ and $w_2(W^4) = w_2(M^6) = w_2 (L_1 \oplus L_2)$.
Hence according to [Thomas1967Z, Lemma1] $W^4$ and $L_1 \oplus L_2$ are  stably isomorphic.  Next  we compute that 
$$e (W^4) = c_ 2(W^4) = c_2 (TM^6, J) ={1 \over 2} ( c_1 ^2 (TM^6, J) - p_1(TM^6) )=  c_1\cdot c_2 = e (L_1 \oplus L_2).$$
Hence,  taking into account [Thomas1967Z, Lemma2]  $W^4$ and $L_1 \oplus L_2$ are isomorphic as real vector  
bundles.  Thus $TM^6$ is the sume of two 3-dimensional vector bundles.\QED

\medskip

{\bf 5.1.2.  Remark. } i)  In 5.3 we discuss regular maximal  non-integrable $G_{\gamma_1}$-structures. If a
$G_{\gamma_1}$-structure is degenerate, but still regular,  then  it is easy to see
that $M^6$ satisfies the condition (P).

ii)  If $M^6$ admits  3 linearly independent vector fields, then it admits also a
$SL(3)\times SL(3)$-structure.  In [Thomas1967I] Thomas give a   necessary and sufficient condition
for an orientable 6-manifold  to admit  3 linearly independent  vector fields, namely $M^6$ has vanishing
Euler class and vanishing Stiefel-Whitney class $w_4$.

\medskip

{\it 5.2. Obstruction for the existence of a $G_{\gamma_2}$-structure. }

\medskip

{\bf 5.2.1. Proposition.} {\it  A manifold $M^6$ admits  a $SL(3,\C)$-structure, if and only if it is orientable
and spinnable.}

\medskip

{\it Proof.}  Clearly 
a 6-manifold $M^6$ admits a $SL(3, \C)$ structure, if and only if  $M^6$ admits an almost complex
structure of vanishing first Chern class.   In particular $M^6$ must be orientable and spinnable.  On the other hand,  if  $M^6$ is orientable and spinnable, then $M^6$ admits $SL(3,\C)$ structure, since 
it admits an almost complex structure, whose first Chern class is an integral lift of $w_2$. Thus the necessary and sufficient condition  for
$M^6$ to admit a $SL(3, \C)$-structure is  the  vanishing of  the  Stiefel-Whitney classes $w_1(M^6)$ and $w_2(M^6)$.\QED

\medskip

{\it 5.3. Maximal non-integrable 3-forms of  $\gamma_1$-type.}  

Every  3-form $\gamma_1$ on $M^6$ defines a pair
two  oriented transversal 3-distributions $D_1$ and $D_2$  together with   volume forms on  each $D_i$  as follows.  Recall that at every point 
$x\in  M$ we can write $\gamma_1 = e^1\wedge e^2\wedge e^3 + e^4\wedge e^5\wedge e^6$.
 The  union
$D_1 \cup D_2$  is   defined uniquely as  the set of all vectors $ v \in T_xM$ such that
$ rank\, (v \rfloor \gamma_1) =2$, or equivalently, $(v\rfloor \gamma_1) ^2 = 0$. The orientation (the volume form)  of $D_1$ and $D_2$ is defined by
the   restriction of $\gamma_1$ to  each distribution $D_i$.  Conversely, a pair of  two  transversal  oriented 3-distributions $D_1$ and
$D_2$ on $M^6$  together with their volume form defines  a  3-form of $\gamma_1$-type  as  follows.  Let their   volume forms be  $\alpha_1$  and $\alpha_2$  respectively. Now we define $\gamma_1 = p_1 ^* (\alpha_1) + p_2^*(\alpha_2)$, where
 $p_1: TM \to D_1$ and $p_2:TM   \to D_2$ are the projections defined by $D_i$.
  
 We   call  structure $(M^6, \gamma_1)$ {\bf regular}, if  the dimensions of  the distributions
$[D_i, D_i] $ defined by $\gamma_1$  are constant over  $M^6$.  We shall call  a regular $G_{\gamma_1}$-structure {\bf maximal non-integrable}, if  at least one of  the distributions
 $D_i$ is maximal non-integrable in the sense that $D_i + [D_i, D_i] = TM$.
 
 At this place we note that the labeling $D_1$ and $D_2$ is well-defined only locally. Globally we may
 be not able to distinguish, what of the two planes is the  $D_1$.  This ambiguity  can be removed, if
 $M^6$ is simply connected, since in this case the two line bundles $\det D_1$ and $\det D_2$
 can be distinguished.
 
 \medskip
 
 We can  describe  the maximal non-integrability of $D_i$ in term of $\gamma_1$ as follows.
 Write $\om_1 = p_1 ^*(\alpha_1), \: \om_2 = p_2 ^* (\alpha_2)$.  Locally we can write
 $\om_1 = p_1^*(e^1\wedge e^2\wedge e^3) $,  $\om_2 = p_2 ^*(e^4\wedge e^5 \wedge e^6)$.

\medskip

{\bf 5.3.1.  Proposition.}  { \it  There is a volume form $D^3 \om_2 \in  \Lambda^3 (\Lambda ^2 (D_1)) ^*$
defined in local  coordinates as follows:
$$ D^3 (\omega_2) = i_1 ^* (d\, p_2^*(e^4)\wedge d\, p_2^*(e^5) \wedge d\,p_2^*( e^6)),$$
where $i_1 : D_i \to TM$ is the embedding, and $dp_2 ^*(e^i)$ are considered as  elements of
$(\Lambda ^2 TM)^*$.
This expression does not depend on the choice of  local 1-forms  $e^i$ considered as 1-forms
on $D_2$.
This volume form  is not zero, if and only if $D_1$ is maximal non-integrable.}

\medskip

{\it Proof.}  We first show that, if $f^4, f^5 , f ^6$ is  another co-frame in $D_2$, so that  $(f^4, f^5, f^6) = g (e^4, e,^5,e^6)$ for $g \in  Gl (D_2)$ then
$$i_1 ^*( d\, p_2^*( f^4)\wedge  d\, p_2 ^*( f^5) \wedge d\,  p_2 ^*(f^6)) = (\det  g) \cdot 
i_1^*( d\,  p_2^*( e^4) \wedge d\,  p_2 ^*( e^5) \wedge d\,  p_2 ^*(e^6)) . \leqno (5.3.2)$$
Proposition 5.3.1  is a local statement, so it suffices to prove on a small  disk  $B^7\subset M^7$.
We denote by $A$  the open  dense subset in the gauge transformation group  $\Gamma(B^7 \times Gl(D_2))$ which is defined by the condition that
$(f^4, e^5, e^6)$ and $(f^4,f^5, e^6)$ are  also  a co-frames on $D_2$.  Then we have
$g = g_3 \circ g_2 \circ g_1$, where $g_1$ sends $(e^4, e^5, e ^6)$ to  $(f^4, e^5, e^6)$,
$g_2$ sends  $(f^4, e^5, e^6)$ to $(f^4, f^5, e^6)$ and $g_3 = g\circ g_1 ^{-1} \circ  g_2 ^{-1}$.
Now it is straightforward to check (5.3.2) for each $g_1$, $g_2$, $g_3$. Hence  (5.3.2) holds
on the  open dense set $A$. Since the  LHS  and RHS of (5.3.2) are continuous mappings, 
the equality (5.3.2) holds  on the  whole  $Gl(D_2)$. This proves the first statement. The second
statement now follows by  direct calculations in local coordinates.\QED

\medskip

  Our study of    maximal non-integrable  $G_{\gamma_1}$-structures is motivated by its   relation with the parabolic  geometry.   This structure  is a generalization  of the famous Cartan 2-distribution in a 5-manifold
and it has  a  canonical conformal structure [Bryant2005].  The Lie algebra of the automorphism group $Aut (M^6, \gamma_1)$
 as well as  local  invariants of $(M^6, \gamma_1)$ can be  calculated explicitly  using the theory of  filtered manifolds (see e.g. [Yamaguchi1993].)

\medskip

\section {Stable 3-forms   on 7-manifolds}

{\it 6.1. Topological  conditions for the existence of a  stable 3-form on a 7-manifold.}

The  sufficient and necessary condition for the existence of   a $G_2$-structure on a 7-manifold $M^7$
has been established by  Gray [Gray1969].  A manifold admits a $G_2$-structure, if and only  if it
is both  orientable and spinnable, i.e. the  first two Stiefel-Whitney classes vanish.

It has been observed in [Le2007]   that a  closed 7-manifold $M^7$  admits a $\tilde G_2$ -structure, if and
only if it is orientable and  spinnable. The closedness condition originates  from
the Dupont work  on obstructions using the K-theory, but  actually what we used  is
 a consequence  of  Dupont's result, namely  the reduction of the $SO(7)$-structure on $M^7$ to a $SO(3)\times SO(4)$-structure.   Thomas in  [Thomas1967M]  also proved that
 any  orientable  closed spin 7-manifold  $M^7$ admits 3-linear independent vector fields. It seems that
 we can drop the condition of closedness  in the proof of his theorem.
 
The geometry of $G_2$-manifolds  has been intensively studied, but the geometry of
$\tilde G_2$-manifolds is barely explored. In [Le2006] we have constructed the first example
of a non-homogeneous   closed 7-manifold which admits a closed 3-form of $\tilde G_2$-type.

\medskip

{\it 6.2.  Malcev algebra structure on  7-manifolds admitting stable 3-forms}.

 Any  stable 3-form $\phi$ in dimension 7 defines  a  reduced  symmetric bilinear form 
by the formula [Bryant1987]
$$<V, W>_\phi = <(V\rfloor \phi)\wedge (W\rfloor\phi)\wedge \phi, \rho>$$
where $\rho$ is some nonzero element in $\Lambda^8(\R^7)$.
 Let us define a multiplication $x\circ_\phi y$ on $\R^7$
 by 
the following  formula:
$$<x\circ_{\phi} y , z>_\phi = \phi (x,y,z).$$
Peter Nagy explained us that the skew-symmetric multiplication $x\circ_\phi y$ defines the structure of the  simple Malcev algebra $A^*$ on $\R^7$ whose corresponding Moufang loop
is $S^7$ for $\phi = \om_1$ in Theorem 4.1 (resp.   the pseudo sphere $S_{(4, 4)}(1)$ of the unit vector
in the  vector space $\R^8$ with  the metric with the signature $(4,4)$ for
$\phi = \om_2)$. Malcev algebras are   generalization of Lie algebras, see [Sagle1961] for more information, in particular the structure of the simple Malcev algebras $A^*$  on $\R^7$.

Thus the tangent bundle $TM^7$ has  the canonical structure of the simple Malcev algebra bundle.

\medskip

\section{Stable 3-forms  on  8-manifolds}

As before we assume that $M^8$ is orientable, since we can go to the orientable double covering, if necessary.

The maximal compact subgroup of $G^+_{\phi_1}$ is $SO(3)$  which is included to $SO(8)$ via the  adjoint representation.  The maximal compact subgroup
of $PSU(1,2)= SU(1,2)/ \Z_3$ is  $S(U(1) \times U(2))/\Z_3$. The subgroups $SO(3)$ and $S(U(1)\times U(2))/\Z_3$ are subgroups of  $PSU(3) = SU(3)/\Z_3$.
Thus any orientable  8-manifold $M^8$ admitting a 3-form of $\phi_1$-type or of $\phi_2$-type  admits also a 3-form of $\phi_3$-type. 
In particular  $M^8$ must be orientable and spinnable. Now  for any spinnable manifold $M^8$ we define the characteristic class  $q_1 (M)$ as follows.

Denote by $q_1$ the spin characteristic class  in $H^4(BSpin (\infty), \Z)$ corresponding to $-c_2 \in H^4(BSU(\infty), \Z)$. For any  spin-bundle $\xi$ over $M$ we  denote by $q_1 (\xi)$ the pull-back
of $q_1$. We set $q_1(M) :=  q_1 (TM)$.

As before  $\rho_2: H^2 (M^8, \Z) \to H^2(M^8, \Z_2)$ denotes
the modulo 2 reduction.  The following Proposition is essentially a reformulation of Corollary 6.4 in
[CCV2007].

\medskip

{\bf 7.1. Proposition}.  {\it  A closed orientable 8-manifold $M^8$ admits a stable 3-form, if and only if it
satisfies the following conditions}
$$ w_2 (M^8) = 0 = e(M^8),\leqno (a)$$
$$  w_6(M^8) \in \rho (H^6 (M^8, \Z)),\leqno(b)$$
$$p_2 (M^8) = - q_1 (M^8) ^2 \text { and }  {( q_1 (M^8) )^2\over 9} [M^8] = 0  \mod 6. \leqno (c)$$

\medskip

 In fact  Corollary 6.4.1 in [CCV2007] is formulated as a necessary and sufficient condition for a manifold
 to admit a $PSU(3)$-structure. But we have seen that the  necessary condition for  a manifold
 $M^8$  to admit a $PSU(3)$-structure is also a necessary condition for a manifold to admit
a $SL(3, \R)$-structure or  a $PSU(1,2)$-structure.

\medskip

\section{Further remarks}

8.1. It is easy to see that  our construction of  natural bilinear forms  works also for   3-forms on space $\R^{3n+2}$. 
In the same way (this  is already noticed first by Bryant for $\R^7$, in [B-V2003]  this  form has been computed for  all except one   multi-symplectic 3-forms) we can associate to any 3-form $\om$ on $\R^{3n +1}$
a bilinear form  with values in $\Lambda ^{3n+1} (\R^{3n+1})^*$, and it descends to a
bilinear  form if the 3-form is non-degenerate; we can also associate to any 3-form
$\om$ on $\R^{3n}$ a linear map from $\R^{3n}$  to $\R^{3n} \otimes \Lambda^{3n} ( \R^{3n}) ^*$,
and this linear map descends to a linear map $\R^{3n} \to \R^{3n}$, if  the 3-form  $\om$ is non-degenerate 
(this is noticed by Hitchin for $\R^6$). We have not yet tested, if  non-degenerate 3-forms exist in higher dimensions. In low dimensions  6,7,8  they  coincide with stable  forms.

\medskip

8.2. Let $\om^3$ be  a stable 3-form   on $M$, $\dim M \ge 7$. Then there is the canonical  inclusion 
$G_\om$ to   $O(k,l)$.  So  if  a  manifold $M$  admits a  stable form $\om^3\not = \gamma_i$, $i =1,2$,  it also
admits a canonical (pseudo)-Riemanian metric.   The curvature of this (pseudo)-Riemannian metric
is a differential invariant of  manifold $(M, \om^3)$. Using these metrics and  existing stable forms we can construct new
differential forms which appear in other special geometries.
Now we shall call a manifold $(M, \om^3)$ stable, if  $\om^3$ is stable.  Stable 
8-manifolds  $(M^8,\om^3)$ seem to us  special interesting, since the bundle $TM^8$  has the canonical  commutative  
multiplication as well as  the structure of
Lie algebra bundle  defined in Proposition 3.11.  
We conjecture that  the algebra $\R^8$ with the commutative  multiplication defined by $\phi_i$ is  a simple algebra. We have a  partial proof for that conjecture  in the case of   $\phi_3$.
The  stable form $\phi_i$ also defines the volume form
on $M^8$ and therefore according to Djokovic it defines   the graded $E_8$-structure on the bundle $\oplus_{i =1}^3(\Lambda^i (TM)\oplus \Lambda^i(T^*M)) \oplus End (TM)$.    

\medskip

8.3. Suppose that $M$ is a compact manifold and $\om^3$ is a stable 3-form on $M$. As we have seen from 8.2  if   $\dim M \ge 7$, then
the automorphism group  $Aut (M, \om^3)$  is a finite dimensional Lie group. If
$\gamma_1$ is maximal non-integrable, then the automorphism group
$(M^6, \gamma_1)$ is alo a finite dimensional Lie group. If $\gamma_1$ is  degenerate, then
the  automorphism group $Aut (M^6, \gamma_1)$ can be infinite dimensional. Example  is
$M^6 = S^1 (\theta^1) \times S^1(\theta^2) \times \Sigma_1\times \Sigma_2$ and
$\om^3 = d\theta^1 \wedge \sigma^1 + d\theta^2\wedge \sigma^2$, where $\sigma^i$ is the volume
element on  the surface $\Sigma_i$. Finally the automorphism group $Aut (M^6, \gamma_2)$ is
also  finite dimensional, since  $SL (3, \C)$ is  elliptic.

\medskip

\section* { Acknowledgement.} We thank Sasha Elashvili  and  Jan Slovak  for    stimulating  discussions
and Peter Nagy for  helpful remarks.
H.-V. L. and J. V.  are partially supported by grant of ASCR Nr IAA100190701.

\medskip

Institute of Mathematics of ASCR, Zitna 25, 11567 Praha 1, Czech Republic, hvle@math.cas.cz\\
Institute of Mathematics of  ASCR, Zizkova 22, 61662 Brno, Czech Republic, naca@ipm.cz\\
Institute of Mathematics of  ASCR, Zizkova 22, 61662 Brno, Czech Republic, vanzura@ipm.cz
\end{document}